\documentclass[12pt]{article}
\usepackage{amsfonts}
\usepackage{latexsym}

\newtheorem{thm}{Theorem}
\newtheorem{cor}{Corollary.}

\newcommand{\PP}{{\mathbb P}}
\newcommand{\qed}{\qquad $ \Box $}

\title{A Remark on the Chisini Conjecture}
\author{Stefan Nemirovski}
\date{}

\begin{document}

\maketitle

In this note, we establish the following consequence of 
Kulikov's results on the Chisini conjecture~\cite{Ku}.

\begin{thm}
\label{main}
A generic ramified covering $f:S\to\PP^2$ of degree at least~$12$
is uniquely determined by its branch curve in\/~$\PP^2.$ 
In other words, the Chisini conjecture holds 
for generic morphisms of degree~$\ge 12$. 
\end{thm}

We discuss the assumptions of this theorem
and the Chisini conjecture itself in section~1. 
Then, in section~2, we present a proof of Theorem~1 
and observe a few related applications of our approach.

\medskip
\noindent
{\bf 1. Generic morphisms of surfaces and the Chisini conjecture.}
A ramified covering $f:S\to\PP^2$ is 
a finite morphism of a non-singular irreducible projective surface 
onto the projective plane.
The {\it branch curve\/} of~$f$ is defined as
the set of points over which $f$ is not \'etale.

A finite morphism $f:S\to\PP^2$ of degree $\deg f\ge 3$
is said to be {\it generic\/} if the following holds:
\begin{itemize}
\item[i)] the branch curve $B\subset\PP^2$ is irreducible 
  and has ordinary cusps and nodes only;
\item[ii)] $f^*B=2R+C$, where the {\it ramification divisor\/}~$R$ 
  is irreducible and non-singular, and $C$ is reduced;
\item[iii)] $f|_R:R\to B$ is the normalization of~$B$.
\end{itemize}

Two generic morphisms $f_1:S_1\to\PP^2$ and $f_2:S_2\to\PP^2$
are called equivalent if there exists an isomorphism~$\varphi:S_1\to S_2$
such that $f_1=f_2\circ\varphi$. 

The precise assertion of Theorem~$1$ is therefore
that two generic morphisms (of {\it a priori\/} different surfaces!) 
having the same branch curve are equivalent 
provided that at least one of them has degree~$\ge 12$.

The above definition of genericity is parallel to the case of Riemann surfaces.
Recall that, according to Riemann and Hurwitz, 
a ramified covering $f:\Sigma\to\PP^1$ is called generic
if over each point in~$\PP^1$ there is at most one quadratic ramification
point of~$f$. 
On the other hand, Theorem~1 shows that the complex
surface case is essentially rigid.

Note that every algebraic surface~$S$ admits a generic morphism $f:S\to\PP^2$
or, in other words, can be represented as a generic ramified covering over~$\PP^2$.
For instance, if $S\subset\PP^r$, then almost every projection
$\PP^r\to\PP^2$ yields a generic morphism~$p:S\to\PP^2$.
Therefore Theorem~1 might be used to understand the moduli
spaces of complex surfaces in terms of the geometry of plane curves.

The conjecture that a generic morphism of degree at least~$5$
is completely determined by its branch curve 
was proposed by O.\,Chisini~\cite{Ch}, who also  
gave an alleged proof of this statement.

The assumption $\deg f\ge 5$ is necessary because of the following
example due to Chisini and Catanese~\cite{Ca}. 
Let $B=C^*$ be the dual curve of a non-singular plane cubic~$C$. 
($B$ is a sextic with $9$ cusps).
Then there exist four non-equivalent generic morphisms 
with branch curve~$B$. Three of them, of degree~$4$, 
are maps from~$\PP^2$ given by projections of the Veronese
surface. 
The fourth map, of degree~$3$, 
is the projection on $\PP^{2*}$ of the elliptic ruled surface obtained 
as the preimage of~$C$ in the incidence variety~$\PP^2\times\PP^{2*}$.
So far this example and its fiber products with ramified coverings 
are the only known examples of non-uniqueness of generic morphisms 
with a given branch curve.

In the important paper~\cite{Mo},
Moishezon proved the Chisini conjecture for branch curves
of generic projections of smooth hypersurfaces in~$\PP^3$
by introducing and analyzing the braid presentations 
of the fundamental group of~$\PP^2\setminus B$.

Recently this approach was substantially developed by Vik.\,Kulikov~\cite{Ku}.
He proved that, for a given branch curve~$B$, the generic morphism
is unique provided that its degree is greater than a certain
function depending on the curve~$B$ (the explicit expression
is given in formula~(\ref{qq}) below).
This allowed him to prove the Chisini conjecture for a wide class
of generic morphism (for instance, for pluri-canonical morphisms 
of surfaces of general type).  

Our (rather modest) contribution is that 
one can estimate Kulikov's function from above 
by using the Bogomolov--Miyaoka--Yau (BMY) inequalities.

In what follows we say that the Chisini conjecture
holds for a class of generic morphisms if every morphism
in this class is determined by its branch curve up to equivalence. 

\medskip
\noindent
{\bf 2. Cusps of branch curves and BMY inequalities.}
Consider a generic morphism $f:S\to\PP^2$ of degree $\deg f=N$
with branch curve $B\subset\PP^2$.
Denote by $2d$ the degree of~$B$ (it is always even), 
by $g$ the genus of the desingularization of~$B$, 
and by $c$ the number of cusps of~$B$.

It was proved in~\cite{Ku} that the morphism $f:S\to\PP^2$
is uniquely determined by~$B$ if
\begin{equation}
\label{qq}
N>\frac{4(3d+g-1)}{2(3d+g-1)-c}.
\end{equation}
We wish to apply the Bogomolov-Miyaoka-Yau inequality on
the algebraic surface~$S$ to estimate the number of cusps 
of the branch curve~$B$ via $g$ and~$d$,
which leads to an {\it a priori\/} upper bound for the right hand side of~(\ref{qq}). To this end we shall need the following formulas 
(cf. Lemmas~6 and~7 in~\cite{Ku}) 
for the self-intersection of the canonical class and the topological Euler characteristic of~$S$:
\begin{equation}
\label{invS}
\begin{array}{rcl}
K_S^2&=&9N-9d+g-1,\\
e(S)&=&3N+2(g-1)-c.
\end{array}
\end{equation}

\smallskip
\noindent
{\it Proof of Theorem\/}~\ref{main}.
Let us assume first that $S$ satisfies the BMY inequality 
\begin{equation}
\label{bmy3}
K_S^2\le 3e(S).
\end{equation}
(This means essentially that $S$ is not a blow-up 
of an irrational ruled surface of irregularity~$>1$.)
Plugging~(\ref{invS}) into the BMY inequality we obtain
$$
9N-9d+g-1\le 9N+6(g-1)-3c,
$$
and therefore
$$
c\le 3d+{\frac53}(g-1).
$$
It follows that
\begin{equation}
\label{chi12}
\frac{4(3d+g-1)}{2(3d+g-1)-c}
\le\frac{12d+4(g-1)}{3d+{\frac13}(g-1)}
=4+\frac{8(g-1)}{9d+(g-1)}<12,
\end{equation}
which proves the Theorem in this case.

Now, if $S$ is (a blow-up of) an irrational ruled surface,
it satisfies the inequality $K_S^2\le 2e(S)$.
(The BMY inequality used above does not follow from this
because both sides may be negative.) 
Arguing in the same way as before, we obtain
$$
c\le -{\frac32}N+{\frac92}d+{\frac32}(g-1) < {\frac32}(3d+g-1).
$$
This gives us the (sharper) estimate
$$
\frac{4(3d+g-1)}{2(3d+g-1)-c}<8,
$$
which completes the proof.\qed

\medskip
In fact, all surfaces of non-general type except~$\PP^2$
satisfy~$K_S^2\le 2e(S)$. 
Moreover, by Theorem~3 in~\cite{Ku} the Chisini conjecture holds 
for generic endomorphisms of~$\PP^2$ of degree~$\ge 5$.
Hence, the second part of our proof yields the following result.

\begin{thm} 
The Chisini conjecture holds for generic morphisms of degree 
at least~$8$ of surfaces of non-general type.
\end{thm}

Returning to the case of surfaces satisfying the ``general"
BMY inequality~(\ref{bmy3}), we see that, in certain cases, 
the last inequality in~(\ref{chi12}) can be improved. 
From the formula for~$K_S^2$ we have
$$
9d=-K_S^2+9N+(g-1).
$$
Plugging this into~(\ref{chi12}) we deduce the following
result (including by the way the case of~$\PP^2$).

\begin{cor}
The Chisini conjecture holds for generic morphisms $f:S\to\PP^2$
of degree~$\ge 8$ such that $K_S^2<9\deg f$.
\end{cor}

Another complementary result can be obtained 
under additional assumptions on the branch curve.

\begin{cor}
The Chisini conjecture holds for generic morphisms of degree~$\ge 8$
such that the genus of the branch curve is less than~$64$.
\end{cor}

\noindent
{\it Proof.}
It follows from the Riemann--Hurwitz formula applied
to the preimage of a projective line that~$\deg f\le d+1$
(cf. Lemma~1 in~\cite{Ku}).
However, if $d\ge 7$ and $g\le 63$ (so that $9d>g-1$), 
then estimate~(\ref{chi12}) improves to~$<8$.\qed

\medskip
\noindent
{\bf Remarks}

\noindent
$1^\circ$ Viktor Kulikov has also obtained the following result:
the Chisini conjecture holds for generic morphisms of degree~$\ge 5$
whose branch curves are cuspidal, i.\,e. have no nodes.

\smallskip
\noindent
$2^\circ$ The estimates for the number of cusps of branch curves
of generic morphisms obtained above are slightly stronger
than the bounds known for arbitrary curves with simple singularities.
It should be noted that the best estimates so far were obtained 
by applying the logarithmic BMY inequality to double coverings
(see~\cite{Hi}).

\medskip
\noindent
{\bf Acknowledgements.}
The author is deeply grateful to Viktor Kulikov for sharing
his results and expertise. This paper was written during the
author's stay at Max-Planck Institute f\"ur Mathematik in Bonn.
It is a pleasure to acknowledge the support of this hospitable institution. 
The author was also supported by RFBR grant 99-01-00969.

\smallskip
{\sc Steklov Mathematical Institute, 117966 Moscow, Russia}

\smallskip
{\it E-mail address:\/} {\tt stefan@mccme.ru}

\end{document}